\documentclass[10pt]{amsart}

\usepackage{graphicx}
\usepackage{amssymb,mathrsfs}

\usepackage{enumerate}        

\usepackage[all,dvips]{xy}    
\UseComputerModernTips        

\newcommand {\comp}[1]{\mathrm{Comp}_{\infty}( #1 )}

\newtheorem{theorem}{Theorem}[section]
\newtheorem{lemma}[theorem]{Lemma}

\newtheorem{proposition}[theorem]{Proposition}
\newtheorem{corollary}[theorem]{Corollary}

\theoremstyle{definition}

\begin{document}

\title[Relative Ends of Thompson's Groups]
{A Proof that Thompson's Groups have Infinitely Many Relative Ends}
\author[D.Farley]{Daniel Farley}
      \address{Miami University of Ohio \\
               Oxford, OH 45056}
      \email{farleyds@muohio.edu}

\begin{abstract}  
We show that each of Thompson's groups $F$, $T$, and $V$ has infinitely
many ends relative to the groups $F_{[0, 1/2]}$, $T_{[0,1/2]}$, and
$V_{[0, 1/2)}$ (respectively).

As an application, we simplify the proof, due to Napier and Ramachandran,
that $F$, $T$, and $V$ are not K\"{a}hler groups. 

We go on to show that Thompson's groups $T$ and $V$ have Serre's property FA.  The main
theorems together answer a question on Bestvina's problem list that was originally posed by 
Mohan Ramachandran.
\end{abstract}



\maketitle

\section{Introduction}

Thompson's group $F$ is the group of piecewise linear homeomorphisms
$h$ of the unit interval such that:  i) each of the finitely many places 
at which $h$ fails to be differentiable are dyadic rational numbers, and
ii) at every other point $x \in [0,1]$, 
$h'(x) \in \{ 2^{i} \mid i \in \mathbb{Z} \}$.  Thompson's groups $T$ and
$V$ have analogous definitions.  The group $T$ is a collection of 
homeomorphisms of the circle, and $V$ can be viewed
as a group of homeomorphisms of the
Cantor set.  A good introduction to all of these groups is \cite{Cannon}.

Recently, Ross Geoghegan posed the problem of determining whether
the group $F$ is K\"{a}hler (see \cite{Brown}).  
A finitely presented group is called
a \emph{K\"{a}hler group} if it is the fundamental group of a compact
K\"{a}hler manifold.  The most important examples of K\"{a}hler groups
(and perhaps the only ones) are the fundamental groups of smooth
complex projective varieties.  

Napier and Ramachandran soon produced proofs
that $F$, $T$, and $V$ are not K\"{a}hler groups \cite{NR2}.  
Their proof in 
\cite{NR1} that
$F$ is not K\"{a}hler used the fact that $F$ is a strictly ascending HNN
extension, and that such groups are never K\"{a}hler.  Their
proofs in \cite{NR2}
that $T$ and $V$ are not K\"{a}hler had two components.  First,
they showed that $T$ and $V$ have infinitely many filtered ends relative
to certain subgroups (see \cite{KR} for the definition of filtered ends).  
Second, they appealed
to the main theorem from \cite{NR2},
which implies that a K\"{a}hler group $G$ having at least $3$ filtered 
ends relative to some subgroup must have a quotient that is isomorphic
to a hyperbolic surface group.  Since the groups $T$ and $V$ are both
simple \cite{Higman}, it is therefore clear that they cannot be K\"{a}hler.

The first half of their argument brought
together a variety of sources, and used the theory of diagram groups
over semigroup presentations \cite{GS}, CAT(0) cubical 
complexes \cite{actions},   
and work of Thomas Klein on filtered ends of pairs of groups \cite{Klein}.  

The main purpose of this note is to simplify the first half of the argument,
and strengthen the conclusion in the process.  We prove the following
theorem (definitions appear in Section \ref{secends}):

\begin{theorem} \label{mrbig}
The pairs $(F, F_{[0, 1/2 ]})$, $(T, T_{[0,1/2]})$, and $(V, V_{[0,1/2)})$
all have infinitely many ends, where $G_{S} = \{ g \in G \mid
g ~is~the~identity~on~ S \}$, for $G \in \{ F, T, V \}$ and $S = [0, 1/2]$
or $[0, 1/2)$. 
\end{theorem}

In Section \ref{secends}, we let $e( G, H)$ denote the number of ends
of the pair $(G,H)$ (or the number of ends of $G$ relative to $H$).  Let
$\tilde{e} (G, H)$ denote the number of filtered ends of the pair $(G,H)$.
The inequality $\tilde{e} (G,H) \geq e( G,H)$ holds true for any group
$G$ and subgroup $H$.  The inequality can be strict \cite{KR}, 
and it is in this 
sense that Theorem \ref{mrbig} strengthens the conclusions of \cite{NR2}. 
Geoghegan \cite{ross} gives examples of pairs for which
$e(G,H) = 3$ but $\tilde{e} (G,H) = \infty$.  Note that we won't  
need to define the 
filtered ends of a group pair here.

Proposition \ref{twoendcase} was originally
proved in \cite{Houghton} (analogous results about filtered ends
were proved in \cite{KR}).  We include our own proof of this Proposition
for the sake of completeness.  As a result, the proof of Theorem \ref{mrbig}
given here is largely 
self-contained, except for the main result of \cite{iso}. 

A second purpose of this note is to show that Thompson's groups $T$ and $V$
both have Serre's property FA, i.e., if $T$ or $V$ acts on a simplicial tree by
automorphisms, then the action has a fixed point.  The proof of this fact in Section \ref{fa}
is due to Ken Brown.  As a consequence, we answer a question posed by Mohan Ramachandran, 
who asked
whether (or to what extent) the following conditions are equivalent for a finitely presented group $G$:  (A)  $G$ has a finite index subgroup
admitting a fixed-point-free action on a simplicial tree, and (B) the pair $(G,H)$ has two or more ends, for some subgroup $H$.  This question
appears on the problem list maintained by Mladen Bestvina.  Our results
show that property  (A) fails for $T$ and $V$, although $T$ and $V$ have multiple (indeed, an infinite number) of ends relative to certain subgroups,
and thus satisfy (B).

I would like to thank Mohan Ramachandran for encouraging me to publish
a proof of Theorem \ref{mrbig}.  The combinatorial approach to group
ends taken in Section \ref{secends} is indebted to \cite{ross}; I thank
Ross Geoghegan for giving me a manuscript version of his book.  After reading
an earlier version of this paper, Mohan Ramachandran told me that Ken Brown had
proved that $T$ and $V$ have property FA, and suggested the relevance of this
fact to the above question.  I thank Ken Brown for his notes (dating from the 1980s), which were the source 
of the argument in Section \ref{fa}.  Portions of this paper were written while I was visiting the Max Planck
Institute for Mathematics.  I thank the Institute for its hospitality and for the excellent working conditions
during my stay. 

\section{Generalities About Ends of Graphs} \label{secends}

Let $\Gamma$ be a locally finite graph, i.e., a locally finite
$1$-dimensional CW complex.  If $C \subseteq \Gamma$ is compact, 
then let $\comp{\Gamma - C}$ denote the set of \emph{unbounded} 
components of 
$\Gamma - C$, i.e., the components having non-compact closure.  
The \emph{number of ends of
$\Gamma$}, denoted $e(\Gamma)$, is
$$ \sup_{C} \{ |\comp{\Gamma - C}| \}.$$

If $G$ is a finitely generated group and $S$ is a finite generating set,
then $\Gamma_{S}(G)$, the \emph{Cayley graph of $G$ with respect to $S$},
is the graph having the group $G$ as its vertex set, and an edge
$e(g,s)$ connecting $g$ to $gs$ for each $g \in G$ and $s \in S$.  The
\emph{coset graph of $H \backslash G$ with respect to $S$}, denoted
$\Gamma_{S}( H \backslash G )$, is the quotient of $\Gamma_{S}(G)$ by the
natural left action of $H$.  

If $G$ is a finitely generated group, then the \emph{number of ends of $G$},
denoted $e(G)$, is the number of ends of its Cayley graph $\Gamma_{S}(G)$,
where $S$ is some finite generating set.  This definition doesn't depend on
the choice of finite generating set, so we will often simply leave off the
subscript $S$, and say that $e(G)$ is the number of ends of $\Gamma(G)$.
In a similar way, we define the \emph{number of ends of the pair $(G,H)$},
denoted $e(G,H)$, by the equation $e(G,H) = e( \Gamma( H \backslash G) )$.

\subsection{A generalization of Hopf's Theorem}

In \cite{Hopf}, Heinz 
Hopf showed that an infinite group has $1$, $2$, or infinitely
many ends.  In this subsection we prove a generalization of this theorem.  
For this it will be helpful to have the following lemma.

\begin{lemma} \label{basic}
If $K$ is a compact subset of the locally finite graph $\Gamma$, then there
is some finite connected subcomplex $K'$ of $\Gamma$ such that
\begin{enumerate}
\item $K \subseteq K'$;
\item $|\comp{\Gamma - K'}| \geq |\comp{\Gamma - K}|$, and
\item each connected component of $\Gamma- K'$ is unbounded. 
\end{enumerate}
\end{lemma}

\begin{proof}
Suppose that $K$ is a compact subset of $\Gamma$.  Let $K_{1}$ be the
the smallest subcomplex of  
$\Gamma$ containing $K$.  It follows from compactness of $K$ that $K_{1}$ is a 
finite subgraph of $\Gamma$.  We enlarge $K_1$ by adding a finite 
number of arcs
to make the resulting graph, $K_2$, connected.   Next, we add all connected
components $C$ of $\Gamma - K_{2}$ having compact closure to $K_2$.  By the
local finiteness of $\Gamma$ and finiteness of $K_2$, the new subgraph
$K_3$ is also compact, and now each component of $\Gamma- K_3$ 
is unbounded.

We set $K_3 = K'$.  It is clear that (1) and (3) are satisfied; we need to
check (2).  Since $K \subseteq K'$, each connected component of 
$\Gamma - K'$ is contained in a (necessarily unique) connected component
of $\Gamma - K$.  If $|\comp{\Gamma - K}| > |\comp{ \Gamma- K'}|$ then there
must be a connected component $C$ of $\Gamma - K$ having non-compact closure
and containing no such connected component $C'$ of $\Gamma - K'$.  Consider
$C - K'$.  The closure $\overline{C - K'}$ is a non-compact, locally finite
graph.  It follows that $C- K'$ contains a connected component $C'$ having
non-compact closure.  Now $C'$ is a connected component of $\Gamma - K'$
and $C' \subseteq C$; this is a contradiction.
\end{proof}

\begin{lemma} \label{basic2}
Let $\Gamma$ be a locally finite graph.  Let $K_1$, $K_2$ be disjoint
finite connected subgraphs such that $\Gamma- K_1$ has $m$ connected
components $C_1$, $\ldots$, $C_m$ and $\Gamma- K_2$ has $n$ connected
components $C'_{1}$, $\ldots$, $C'_{n}$.  If $K_1 \subseteq C'_{1}$
and $K_2 \subseteq C_1$, then $C_{2}$, $\ldots$, $C_{m}$, $C'_{2}$,
$\ldots$, $C'_{n}$ are distinct connected components of
$\Gamma - ( K_{1} \cup K_{2} )$.
\end{lemma}

\begin{proof}
We first show that $C_{2}, C_{3}, \ldots, C_{m}, C'_{2}, C'_{3},
\ldots, C'_{n}$ are in fact components of $\Gamma - ( K_{1} \cup K_{2})$.
Choose a component $C_{i}$ of $\Gamma - K_{1}$  $(2 \leq i \leq m)$.
Pick an arbitrary point $x \in C_{i}$;  
let $\widehat{C}$ be the component of $\Gamma - (K_{1} \cup K_{2})$
containing $x$.
Pick an arbitrary point $y \in C_{i}$.  Let $p_{xy} \subseteq C_{i}$
be the image of a path connecting $x$ to $y$.  Since $p_{xy} \subseteq C_{i}$,
$p_{xy} \cap K_{1} = \emptyset$ (because $C_i \subseteq \Gamma - K_{1}$)
and $p_{xy} \cap K_{2} = \emptyset$ (because $K_{2} \subseteq C_{1}$ and
$C_i \cap C_1 = \emptyset$).  It follows that 
$C_{i} \subseteq \widehat{C}$.
To prove the reverse inclusion, 
let $y \in \widehat{C}$ be arbitrary.  Let $p_{xy} \subseteq 
\widehat{C}$ be the image of a path connecting $x$ to $y$.  It follows
that $p_{xy} \subseteq \Gamma - ( K_{1} \cup K_{2} ) \subseteq 
\Gamma - K_{1}$.  This directly implies that $y \in C_{i}$, so
$\widehat{C} = C_{i}$.
It follows (by symmetry in the case of $C'_{2}, \ldots, C'_{n}$)
that
each of $C_{2}, \ldots, C_{m}, C'_{2}, \ldots, C'_{n}$ is a connected 
component
of $\Gamma - ( K_{1} \cup K_{2})$.

It is easy to see that $C_{2}, \ldots, C_{m}$ are distinct components
of $\Gamma - (K_{1} \cup K_{2})$, and that $C'_{2}, \ldots, C'_{n}$
are also distinct components of $\Gamma - (K_{1} \cup K_{2})$.  
Suppose $C_{i} = C'_{j}$ ( for some $2 \leq i \leq m$, 
and $2 \leq j \leq n$).  Let
$x \in C_{i}$ and let $y \in C_{i'}$, where $i' \neq i$ and $2 \leq i' \leq m$.
By the connectedness of $K_1$, there is a path whose image $p_{xy}$
satisfies $p_{xy} \subseteq K_{1} \cup C_{i} \cup C_{i'}$.  It follows
that $p_{xy} \subseteq \Gamma - K_{2}$, so in fact $p_{xy} \subseteq
C'_{j}$.  But now it follows that $p_{xy} \subseteq \Gamma - K_{1}$
(since $C'_{j} \cap K_{1} \subseteq C'_{j} \cap C'_{1} = \emptyset$),
a contradiction.  
\end{proof}

\begin{theorem} \label{hopf} If $\Gamma$ is a locally finite graph admitting an
infinite group of covering transformations, then $e( \Gamma) = 1, 2,$ or
$\infty$.
\end{theorem}

\begin{proof}
Suppose that $3 \leq e( \Gamma) < \infty$.  It follows that there is some
compact set $K \subseteq \Gamma$ such that
 $|\comp{\Gamma - K}| = e( \Gamma) =n$.
By the previous Lemma \ref{basic}, 
we may assume that $K$ is a finite connected subcomplex
of $K$, and that each component of $\Gamma - K$ is unbounded.

Since $\Gamma$ admits an infinite group of covering transformations, there
must exist some covering transformation $\gamma$ such that
$\gamma \cdot K \cap K = \emptyset$.  Let $C_{1}, \ldots, C_{n}$ denote
the connected components of $\Gamma - K$; let $C'_{1}, \ldots, C'_{n}$
denote the connected components of $\Gamma - \left( \gamma \cdot K \right)$.
All of these connected components are unbounded by our assumptions.
By the connectedness of $K$, we can assume, without loss of generality,
that $K \subseteq C'_{1}$.  Similarly $\gamma \cdot K \subseteq C_{1}$,
without loss of generality.  It follows that $C_{2}, \ldots, C_{n}, C'_{2},
\ldots, C'_{n}$ are $2n-2$ unbounded components of 
$\Gamma - \left( K \cup \gamma \cdot K \right)$.  
Since $3 \leq e( \Gamma) = n$ and $2n - 2 \leq e( \Gamma)$, we have a
contradiction.
\end{proof}

\begin{corollary} \label{onetwoinf}
If $G$ is a finitely generated group, then $e(G) = 0, 1, 2,$ or $\infty$.
If $H \leq G$ has infinite index in its normalizer 
$N_{G}(H) = \{ g \in G \mid gHg^{-1} = H \}$, then $e(G,H) = 1, 2,$ 
or $\infty$.
\end{corollary}

\begin{proof}
The first statement follows easily after applying the previous theorem
to the Cayley graph $\Gamma(G)$.  (The case $e(G)=0$ corresponds to the
case in which $G$ is finite.)  The second statement follows from applying
Theorem \ref{hopf} 
to the coset graph $\Gamma( H \backslash G)$ and noticing that
$N_{G}(H)/ H$ acts as covering transformations on $\Gamma( H \backslash G)$.
\end{proof}

\subsection{The set of ends of a graph}

In certain situations, it is useful to work with a set of ends, rather
than simply a number of ends.
Let $c: [0, \infty) \rightarrow \Gamma$ be a cellular proper ray, i.e.,
each open interval $(i, i+1)$ (for $i \in \mathbb{Z}$) 
is mapped homeomorphically to an open edge by
$c$, and $c^{-1}(K')$ is compact if $K'$ is.  Two cellular proper rays
$c$ and $c'$ are joined by a \emph{proper ladder} if there is
a map 
$$L: \left( [0, \infty) \times \{ 0 \} \right) \cup 
\left( [0, \infty) \times \{ 1 \} \right) \cup \left(Z     
\times [0,1] \right) \rightarrow \Gamma$$
where $Z$ is an infinite subset of the positive integers,
$L(t,0) = c(t)$, $L(t,1) = c'(t)$, and $L^{-1}(K')$ is
compact if $K'$ is.  We say that $c$ and $c'$ \emph{define the same end},
and write $c \sim c'$, if $c$ and $c'$ are joined by a proper ladder.  It
is rather clear that $\sim$ is an equivalence relation; the equivalence 
classes are called \emph{ends}.  The set of ends of $\Gamma$ is denoted 
$\mathcal{E}( \Gamma)$.  

\begin{proposition} \label{ends}
Let $c_{1}$, $c_{2}$ be proper rays in the locally finite graph $\Gamma$.
The proper rays $c_{1}$, $c_{2}$ define the same end if and only if for
any compact subset $K$ of $\Gamma$ there is some $t \in \mathbb{R}$ so
that $c_{1} ( [ t, \infty))$ and $c_{2} ( [t, \infty))$ are in the
same component of $\Gamma - K$.
\end{proposition}

\begin{proof}
($\Rightarrow$) Suppose that $c_1$ and $c_2$ define the same end, and let
$K$ be a compact subset of $\Gamma$.  Let $\ell$ be  a proper ladder
joining $c_1$ and $c_2$.  Choose $t_1$ large enough that
$c_{1} ( [t_{1}, \infty )) \cap K = c_{2} ( [ t_{1}, \infty )) \cap K
= \emptyset$.  The properness of $\ell$ implies that some rung 
$\ell ( \{ t \} \times [0,1] )$ of the ladder for $t > t_1$ is disjoint
from $K$, and this directly implies that
$c_{1} ( [ t, \infty ))$ and $c_{2} ( [t , \infty))$ are in the same
component of $\Gamma - K$.  

($\Leftarrow$) We need to connect $c_1$ to $c_2$ by a proper ladder.
Begin by connecting $c_{1}(1)$ to $c_{2}(1)$ by an arbitrary arc, to form
the first rung $\ell ( \{ 1 \} \times [0,1] )$ of the ladder.  Now 
we can use the hypothesis with $K = c_{1} ( [0,1] ) \cup c_{2} ( [0,1] )
\cup \ell ( \{ 1 \} \times [0,1] )$ to add another rung to the ladder, 
which is disjoint from the first rung.  By continuing in the same way, we
inductively define a proper ladder between $c_1$ and $c_2$.
\end{proof}

\begin{corollary}
If $\Gamma$ is a locally finite graph and $e( \Gamma ) = m >0$, then
there is some compact subset $K$ of $\Gamma$ such that $\Gamma - K$ has
exactly $m$ connected components $C_{1}, \ldots, C_{m}$, all of which are
unbounded.
Let $K$ be any such compact subset. Two proper rays $c_1$, $c_2$
represent the same end if and only if $c_{1} ([ t, \infty)), 
c_{2}( [ t, \infty )) \subseteq C_{i}$ for some $i$ and sufficiently 
large $t$.

In particular, if $e( \Gamma)$ is finite, then 
$e( \Gamma ) = | \mathcal{E} ( \Gamma ) |$.
\end{corollary}

\begin{proof}
The existence of $K$ is an immediate consequence of the statement that
$e( \Gamma ) = m$ and Lemma \ref{basic}.  
The forward direction of the second statement is clear.

Suppose $c_1$ and $c_2$ are two proper rays and
$c_{1} ( [ t, \infty )), c_{2} ( [ t, \infty )) \subseteq C_i$
for some $i$ and some $t$.  If $c_1$ and $c_2$ define
separate ends, then by Proposition \ref{ends} 
there is some compact $K' \subseteq \Gamma$ so
that, for some $t'$, $c_{1} ( [t', \infty ))$ and $c_{2} ( [ t', \infty))$
are contained in distinct components of $\Gamma - K'$.  It follows from this 
that $C_{i}$ contains two unbounded components $C'$, $C''$ of $\Gamma - K'$.
This implies that 
$|Comp_{\infty} ( \Gamma - K' )| >  
|Comp_{\infty} ( \Gamma - K )| = e(\Gamma)$, a contradiction.
\end{proof}   

\subsection{The case of two ends}

\begin{proposition} \label{twoendcase}
Let $\Gamma$ be a locally finite graph admitting an infinite group 
$C(\Gamma)$ of covering transformations.  If $e(\Gamma) =2$, then
$C(\Gamma)$ has an infinite cyclic subgroup of finite index.
\end{proposition}

\begin{proof}
Suppose $e(\Gamma) = 2$.  Let $K$ be a finite connected subgraph of
$\Gamma$ such that $\Gamma - K = C_{1} \cup C_{2}$, where each
$C_{i}$ is an unbounded connected component of $\Gamma - K$.  The
group $C ( \Gamma )$ acts on the set of ends, and after passing to
a subgroup of index $2$ if necessary, we can assume that $C(\Gamma)$
fixes both ends.  Since $\Gamma$ is infinite, locally finite, and $C(\Gamma)$
acts freely, there is some $\gamma \in C( \Gamma)$ so that
$(\gamma \cdot K) \cap K = \emptyset$.  Assume without loss of generality
that $(\gamma \cdot K) \subseteq C_{1}$.

We first show that $(\gamma \cdot C_{1}) 
\cap C_{2} = \emptyset$.  Suppose $x \in (\gamma \cdot C_{1}) \cap C_{2}$;
let $y \in (\gamma \cdot C_{2}) \cap C_{2}$ 
(Here 
$ (\gamma \cdot C_2) \cap C_2 \neq \emptyset$ since $\gamma$ fixes the ends
of $\Gamma$).  Since $x$ and $y$ are in $C_{2}$, there is an edge-path
$p$ connecting $x$ to $y$ in $C_{2}$, and
$p \cap (\gamma \cdot K) \subseteq p \cap C_{1} = \emptyset$.  Thus
$(\gamma^{-1} \cdot p)$ is an edge-path $p$ connecting 
$(\gamma^{-1} \cdot x)
\in C_{1}$ to $(\gamma^{-1} \cdot y) \in C_{2}$ and missing $K$.  This
is a contradiction, so $( \gamma \cdot C_{1} ) \cap C_{2} = \emptyset$.

Now note that $K \subseteq \gamma \cdot C_{2}$.  For otherwise
$K \subseteq \gamma \cdot C_{1}$, and since $\gamma \cdot C_{1}$
would then be an open set containing $K$, it would follow that 
$\gamma \cdot C_{1}$ contains elements $C_{2}$, a contradiction.

We can now apply Lemma \ref{basic2}:  
since $(\gamma \cdot K) \subseteq C_{1}$ and
$K \subseteq (\gamma \cdot C_{2})$, $C_{2}$ and $(\gamma \cdot C_{1})$ are 
distinct connected components of $\Gamma - ( K \cup (\gamma \cdot K ))$,
and both are clearly unbounded.  

We have
\begin{eqnarray*}
       \Gamma - ( K \cup (\gamma \cdot K) )  & = &
( C_{1} \cap (\gamma \cdot C_{1}) ) \cup ( C_{2} \cap (\gamma \cdot C_{1}) )
\cup ( C_{1} \cap (\gamma \cdot C_{2}) ) \\ 
& & \cup ( C_{2} \cap (\gamma \cdot C_{2}) ) \\
& =  & 
C_{2} \cup (\gamma \cdot C_{1}) \cup ( C_{1} \cap (\gamma \cdot C_{2}) ).\\
\end{eqnarray*}
Indeed, the first equality above is obvious.  The forward inclusion of  
the second equality follows from the fact that 
$( \gamma \cdot C_{1}) \cap C_{2} = \emptyset$.  The reverse inclusion
follows from Lemma \ref{basic2}: since $\gamma \cdot K \subseteq C_{1}$
, $K \subseteq ( \gamma \cdot C_{2} )$, and $K$ is compact and connected,
$C_{2}$ and $( \gamma \cdot C_{1} )$ are distinct connected components
of $\Gamma - ( K \cup \gamma \cdot K )$ by the Lemma.  Thus 
$C_{2} \cup ( \gamma \cdot C_{1} ) \subseteq 
 \Gamma  - ( K \cup \gamma \cdot K )$, and the reverse inclusion is 
established.   

The second equality above easily implies that
$C_{2} \subseteq  \gamma \cdot C_{2}$ and 
$\gamma \cdot C_{1} \subseteq C_{1}$.  Indeed, both of these last 
inclusions are proper:  the first is proper since 
$K \subseteq \gamma \cdot C_{2}$, and the second is proper since
$\gamma \cdot K \subseteq C_{1}$.  This directly implies that
$\gamma$ has infinite order. 
Note also that $\overline{ ( C_{1} \cap (\gamma \cdot C_{2} ))}$
is compact, since $e ( \Gamma ) = 2$ and $C_{2}$, $\gamma \cdot C_{1}$
are both unbounded.

Next we need to show that, for any $x \in C_{1}$, there is $n <0$ such 
that $\gamma^{n} \cdot x \in C_{2}$, and  for any $x \in C_{2}$,
there is $n > 0$ such that $\gamma^{n} \cdot x \in C_{1}$.  We argue by
contradiction:  suppose $x \in C_{1}$ and
$\{ \gamma^{-1} \cdot x , \ldots , \gamma^{-n} \cdot x , \ldots \} \subseteq 
K \cup C_{1}$.  Let $y \in C_{2}$; choose some path $p$ connecting $x$ to $y$.
Now $\{ \gamma^{-1} \cdot y , \ldots , \gamma^{-n} \cdot y , \ldots \}
\subseteq C_{2}$, so each path $\gamma^{-1} \cdot p , \ldots , \gamma^{-n}
\cdot p, \ldots$ meets $K$.  It follows that some subsequence of  
$\gamma^{-1} \cdot x , \ldots , \gamma^{-n} \cdot x , \ldots $ has a limit
(by the local finiteness of $\Gamma$), and thus infinitely
many terms of the sequence are identical (since
$\gamma$ is a covering transformation), which implies that 
$\gamma^{k} = 1$ for some $k \neq 0$.  This is a contradiction.  Thus, 
for any $x \in C_{1}$, there is $n < 0$ so that $\gamma^{n} \cdot x \in
C_{2}$.  A similar argument shows that, for any $x \in C_{2}$, there is
$n > 0$ so that $\gamma^{n} \cdot x \in C_{1}$.

Finally, we argue that
$\widehat{K} = K \cup ( \gamma \cdot K) \cup 
( C_{1} \cap (\gamma \cdot C_{2}))$ is a compact fundamental domain for the
action of $\langle \gamma \rangle$ on $\Gamma$; a standard argument then
shows that $C( \Gamma)$ contains $\langle \gamma \rangle$ as a finite-index
subgroup.  Compactness of $\widehat{K}$ has already been established.  
Let $x \in \Gamma - K$.  We may assume that $\left( \langle \gamma \rangle
\cdot x \right) \cap K = \emptyset$.  The argument of the previous paragraph
shows that there is some $n_{1}< 0$ so that $\gamma^{n_{1}} \cdot x \in C_{2}$
and some $n_{2} > 0 $ so that $\gamma^{n_{2}} \cdot x \in C_{1}$.  
Thus, there are consecutive integers $k$, $k+1$
such that $\gamma^{k} \cdot x \in C_{2}$ and $\gamma^{k+1} \cdot x \in C_{1}$.
It follows that $\gamma^{k+1} \cdot x \in C_{1} \cap ( \gamma \cdot C_{2})$.
\end{proof}

\subsection{Almost Invariant Subsets}

In the main argument of this paper, we will need a criterion, due to
Sageev, for the pair $(G,H)$ to have multiple ends.  If $G$ acts on a
set $S$, then a subset $T$ of $S$ is said to be \emph{almost invariant}
if the symmetric difference 
$|gT \triangle T|$ is finite for any $g \in G$. 

\begin{theorem} \label{sageev} \cite{Sageev} 
Let $G$ be a finitely generated group; let $H \leq G$.  
Consider the left (or right) 
action of $G$ on the set $G/H$ of left (or right) cosets of $H$.
If there is a subset $A$ of $G/H$ such that
\begin{enumerate}
\item $A$ is almost invariant, and
\item each of $A$, $A^{c}$ is infinite,
\end{enumerate}
then $e(G,H) \geq 2$.  Conversely, given a pair $(G,H)$ such that 
$e(G,H) \geq 2$, there exists such an almost invariant set $A$ in $G/H$.
\end{theorem}

\begin{proof}
($\Rightarrow$) 
We prove the theorem in the case of the right action of $G$ on the
collection of right cosets $H \backslash G$.  Begin by choosing a 
finite generating set $S$ 
for $G$, and building the coset graph
$\Gamma_{S}( H \backslash G)$.  Suppose there is a set $A$ as in the
statement of the theorem.

The elements of $A$ can naturally be identified with vertices of
$\Gamma_{S} ( H \backslash G)$.  Consider the collection of all edges
$e$ which connect an element of $A$ with an element in $A^{c}$.  We
claim that the union $K$ of all such edges is a finite subgraph
of $\Gamma_{S} ( H \backslash G)$.  If not, then, by the finiteness of
$S$ and without loss of generality, there must be infinitely
many disjoint directed
 edges of $K$, each labelled by the same generator $s \in S$, and
each running from an element of $A^{c}$ to an element of $A$.  It
follows from this that each of the (infinitely many) terminal vertices of
such edges are in $A$, but not in $As$.  This implies that $A \triangle
As$ is an infinite set, which contradicts the fact that $A$ is almost
invariant.

Since $A - K$ and $A^{c} - K$ are both infinite sets, and they are clearly
separated by $K$, it follows that $\Gamma_{S} ( G, H)$ has at least
two ends. 

($\Leftarrow$)  We won't need to use this implication, so we leave
the (easy) proof as an exercise.  The idea is to choose  a compact subgraph
$K$ which divides $\Gamma$ into at least two unbounded components, and then
use the vertices of one of these components as $A$.
\end{proof}

\section{A proof that Thompson's groups have infinitely many relative ends}  

\begin{lemma} \label{trivial}
Let $G_{1}$ and $G_{2}$ be finitely generated groups.  Suppose that
$G_{1} \leq G_{2}$, $H_{1} \leq G_{1}$, $H_{1} \leq H_{2}$, 
and $H_{2} \leq G_{2}$.    

If the natural map $\phi: G_{1} / H_{1} \rightarrow G_{2} / H_{2}$ 
is injective
and $A \subseteq G_{2}/H_{2}$ is an almost invariant subset (under the left
action of $G_{2}$), then $\phi^{-1} (A)$ is almost invariant under the left
action of $G_{1}$.
\end{lemma}

\begin{proof}
Let $A$ be an almost invariant subset of $G_{2} / H_{2}$ under the left
action of $G_{2}$.  We consider the inverse image $\phi^{-1}(A)$; let
$g \in G_{1}$.  We have that 
$$ \phi ( \phi^{-1}(A) \triangle g \phi^{-1}(A) ) 
\subseteq  A \triangle g A.$$
Since $\phi$ is injective, it directly follows that $\phi^{-1}(A) \triangle
g \phi^{-1}(A)$ is finite.
\end{proof}

\begin{proposition} \label{big}
Let $V_{[0,1/2)}$ denote the subgroup of Thompson's group $V$ which acts
as the identity on $[0,1/2)$.  
\begin{enumerate}
\item The set $A = \{ g V_{[0,1/2)} \mid g_{\mid_{[0,1/2)}} ~is~affine~
\}$ is almost invariant under the action of $V$ on 
$V/ V_{[0,1/2)}$.  Both $A$ and its complement are infinite.
\item The quotient group 
$N\left( V_{[0,1/2)} \right) /  
 V_{[0,1/2)}$ 
has no cyclic subgroup of finite
index.
\end{enumerate}
In particular, $e( V, V_{[0,1/2)}) = \infty$.
\end{proposition}
  
\begin{proof}
(1)  The statement that $A$ is almost invariant is essentially the content
of \cite{iso}.  
The main argument of \cite{iso} shows that $(v-1) \cdot \chi_{A}$ (where
$\chi: P(V / V_{[0,1/2)}) \rightarrow \mathbb{Z}$ 
is the characteristic function)
is a finite
sum for any element $v \in V$.  This clearly means that $A$ is almost 
invariant. 

For suitable selections
of elements $g_i$ ($i$ a positive integer), $g_{i}$ is affine on $[0,1/2)$ and
$g_{i} \cdot [0,1/2)$ is  
the dyadic interval $[0, 2^{-i})$.  For instance, we can let
$g_{i} = x_{0}^{i}$, where $x_{0}$ is one of the standard generators
of $F \subseteq V$.  As a piecewise linear homeomorphism of $[0,1]$, 
$x_0$ is 
defined as follows:
$$ x_{0}(t) = \Bigg\{ \begin{array}{ll} \frac{1}{2} t & 0 \leq t \leq 1/2 \\
t - \frac{1}{4} & 1/2 \leq t \leq 3/4 \\
2t-1 & 3/4 \leq t \leq 1
\end{array}$$
\noindent The cosets $g_{i} V_{[0, 1/2)}$ are
easily seen to be distinct, so $A$ is infinite.

Any two distinct 
elements of the infinite subgroup $V_{[1/2, 1)}$ represent distinct 
left cosets of $V_{[0,1/2)}$, and only one of these left cosets (containing
the identity) lies in $A$.  It follows that $A^{c}$ is infinite.  
This proves (1). 

(2)  
Each element of $V_{[1/2, 1)}$ normalizes $V_{[0,1/2)}$, and any two elements
in $V_{[1/2, 1)}$ represent distinct left cosets of $V_{[0,1/2)}$.  It follows
that the group $V_{[1/2, 1)}$ embeds in the quotient from the statement of
the proposition.  But $V_{[1/2, 1)}$ is isomorphic to $V$ itself, and $V$
has no cyclic subgroup of finite index.  This proves (2).

The final statement now follows from (1), (2), Theorem \ref{sageev},
Corollary \ref{onetwoinf}, and Proposition \ref{twoendcase}. 
\end{proof}

\begin{proposition} 
$e(T, T_{[0,1/2]}) = e(F, F_{[0,1/2]}) = \infty$.
\end{proposition}

\begin{proof}
We argue that $e(F, F_{[0,1/2]}) = \infty$.  The case of the group $T$
is similar.

We first check the hypotheses of Lemma \ref{trivial}.  It is clear
that there are inclusions $F \rightarrow V$, $F_{[0,1/2]} \rightarrow 
V_{[0,1/2)}$.  We next have to show that the induced map
$$ \phi: F / F_{[0,1/2]} \rightarrow V / V_{[0,1/2)}$$
is injective.  Suppose that $g_{1} F_{[0,1/2]}$
and $g_{2} F_{[0,1/2]}$ both have the same image under $\phi$.  It follows
that $g_{2}g_{1}^{-1} \in V_{[0,1/2)}$.  Now clearly $g_{2}g_{1}^{-1} \in F$,
and it then follows from continuity that $g_{2}g_{1}^{-1} \in F_{[0,1/2]}$,
so $g_{1} F_{[0,1/2]} = g_{2} F_{[0,1/2]}$.  Therefore $\phi$ is injective.
This implies that 
$\phi^{-1}(A) = \{ g F_{[0,1/2]} \mid g_{|_{[0,1/2]}}~is~linear \}$
is an almost invariant set.  
 
We can prove that $\phi^{-1}(A)$ and $\phi^{-1}(A)^{c}$ are both infinite sets 
as in the proof of Proposition \ref{big}.  Indeed, as before,
the cosets $x^{i}_{0} F_{[0,1/2]}$ are all distinct (proving that
$\phi^{-1}(A)$ is infinite) and any two distinct elements
of $F_{[1/2, 1]}$ define distinct cosets of $F_{[0,1/2]}$, and exactly
one of these cosets ($F_{[0,1/2]}$ itself) is in $\phi^{-1}(A)$.  This
proves that $\phi^{-1}(A)^{c}$ is also an infinite set.  It follows
that $e( F, F_{[0,1/2]}) = 2$ or $\infty$.

As in Proposition \ref{big}, there is an embedding of $F_{[1/2, 1]}$
into the quotient group 
$$ N( F_{[0,1/2]}) / F_{[0,1/2]},$$
and $F_{[1/2, 1]} \cong F$.  Since $F$ has no infinite cyclic subgroup of
finite index, it follows that $e( F, F_{[0,1/2]}) = \infty$.
\end{proof} 

\section{Proof that $T$ and $V$ have Serre's property $FA$} \label{fa}

Suppose that $G$ acts simplicially on the simplicial tree $\Gamma$.  We say
that $G$ acts \emph{without inversions} if, whenever $g \in G$ leaves an edge $e$ invariant, 
$g$ acts as the identity on $e$.  We will assume (after barycentrically subdividing the tree, if necessary)
that any simplicial action of $G$ on a tree is an action without inversions.  
We say that $G$ has \emph{property $FA$} if every simplicial action     
of $G$ on a tree has a fixed point, i.e., $G\cdot v = v$, for some $v \in \Gamma$. 

We will need some standard facts about automorphisms of trees.  
If $g \in G$, then 
$Fix( g) = \{ x \in \Gamma \mid g\cdot x = x \}$.  The set $Fix(g)$ is a subtree
of $\Gamma$ if it is non-empty.  If $Fix(g) \neq \emptyset$,
then $g$ is called \emph{elliptic}; otherwise, $g$ is \emph{hyperbolic}.

\begin{lemma} \label{trees}
Let $G$ be a group acting on a simplicial tree $\Gamma$ by automorphisms.
\begin{enumerate}
\item Let $g \in G$.  Either $g$ acts on a unique simplicial line in $\Gamma$ 
by 
translation  
(called an \emph{axis} for $g$), 
or $Fix(g) \neq \emptyset$.
\item If the fixed sets $Fix(g_{1}), Fix(g_{2})$ are non-empty and disjoint, then $Fix(g_{1}g_{2}) = \emptyset$.
\item If $g_{1}$ and $g_{2}$ are elliptic and $g_{1} \cdot Fix(g_{2}) = Fix(g_{2})$, then $Fix(g_{1}) \cap Fix(g_{2}) \neq \emptyset$.
\item If $G$ is generated by a finite set of elements $s_{1}, \ldots, s_{m}$ such that
the $s_{j}$ and the $s_{i}s_{j}$ have fixed points, then $G$ has a fixed point.
\end{enumerate}
\end{lemma}

\begin{proof}
(1) is a consequence of Proposition 24 (page 63) from \cite{Serre}.  Our proof of (2) uses Corollary 1 from page 64 of 
\cite{Serre}, which says that if $abc = 1$ and each of $a, b, c$ is elliptic, then $a, b,$ and $c$ have a common fixed point. 
If we read this corollary with $g_{1} =a$, $g_{2} = b$, and $g_{2}^{-1}g_{1}^{-1} = c$, then the assumption that $g_{1}g_{2}$
is elliptic leads to a contradiction, since $Fix(g_{1}) \cap Fix(g_{2}) = \emptyset$.  To prove (3), we use Lemma 9 from page 61 of
\cite{Serre}, which asserts that, if $\Gamma_{1}$ and $\Gamma_{2}$ are disjoint subtrees of $\Gamma$, then there is a unique minimal geodesic
segment $\ell$ connecting $\Gamma_{1}$ to $\Gamma_{2}$.  That is, if $\ell'$ connects a vertex of $\Gamma_{1}$ with a vertex in $\Gamma_{2}$,
then $\ell \subseteq \ell'$.  It is fairly clear that $\ell$ meets each of the subtrees $\Gamma_{1}$ and $\Gamma_{2}$ in exactly one point.
Now suppose that $g_{1}$ and $g_{2}$ are elliptic; we assume that their fixed sets are disjoint. 
Let $\ell = [x,y]$ connect $x \in Fix(g_{1})$ with $y \in Fix(g_{2})$.  We assume that $\ell$ is the minimal geodesic
connecting these fixed sets.  Our assumptions
imply that $[x,y] \cup g_{1} \cdot [x,y]$ is a geodesic segment connecting two points in $Fix(g_{2})$, namely $y$ and $g_{1} \cdot y$.  Since
$Fix(g_{2})$ is a tree, $[x,y] \cup g_{1} \cdot [x,y] \subseteq Fix(g_{2})$.  This implies the contradiction $Fix(g_{1}) \cap Fix(g_{2})
\neq \emptyset$, proving (3).  

Statement (4) is Corollary 2 from \cite{Serre}, page 64.
\end{proof}
    
Throughout the rest of this section, we assume that $G = T$ or $V$.  Let $g \in G$.
We say that $g$ is \emph{small} if $g$ is the identity on some
\emph{standard dyadic subinterval of $[0,1]$} (i.e., some subinterval of the form $[ \frac{i}{2^{n}} , \frac{i+1}{2^{n}} ]$, where
$n$ is a non-negative integer and $0 \leq i < 2^{n}$). 

\begin{lemma}  
If $G$ acts on a tree $\Gamma$, and $g \in G$ is small, then $g$ has a fixed point.
\end{lemma}

\begin{proof}
Let $g \in G$.  Suppose, for a contradiction, that $g$ is hyperbolic and acts by translation on the geodesic line $\ell$.  
Let $I$ be a standard dyadic subinterval of $[0,1]$ such that $g_{\mid I} = id_{I}$.  We consider the subgroup
$H$ of $G$ having support in $I$.  (This group is isomorphic either to $F$ or to $V$, depending on whether
$G$ is $T$ or $V$, respectively.)  For any $h \in H$, $hgh^{-1} = g$, so
$$  g \cdot h \ell = hgh^{-1} \cdot h \ell = h \cdot g \ell = h \cdot \ell.$$
It follows from the uniqueness of the axis $\ell$ that 
$h \ell = \ell$.  Thus, the entire group $H$ leaves the line $\ell$ invariant, so there
is a homomorphism $\phi: H \rightarrow D_{\infty}$.  The kernel of $\phi$ is large:  it will
contain $[H,H]$ (if $G =T$, in which case $H$ is isomorphic to $F$, and every proper quotient of $F$ is abelian by Theorem
4.3 from \cite{Cannon}) or $H$ itself (if $G = V$, in 
which case $H$ is isomorphic to $V$, and therefore simple).  By Theorem 4.1 from \cite{Cannon}, $[H,H]$ consists
of the subgroup of $H$ which acts as the identity in neighborhoods of the left- and right-endpoints of $I$.
Thus, if $I'$ is a standard dyadic interval whose closure is contained in the interior $I$, then 
any element of $G$ supported in $I'$ will fix the entire line $\ell$.  
  
We can choose a conjugate $kgk^{-1}$, where $k \in G$, and the support of $kgk^{-1}$ is contained in $I'$.  It follows
that $kgk^{-1}$ is elliptic, and therefore $g$ is elliptic as well.  This is a contradiction. 
\end{proof}

\begin{theorem}  
$T$ has Serre's property FA.
\end{theorem}

\begin{proof}
Identify $[0,1]/ \sim$ with the standard unit circle $S^{1} = \{ (x,y) \in \mathbb{R}^{2} \mid x^{2} + y^{2} = 1 \}$ by
the quotient map $f: [0,1] \rightarrow S^{1}$, where $f(t) = ( \cos(2\pi t), \sin(2\pi t) )$.  This identification induces
an action of $T$ on $S^{1}$.  We consider four subgroups of $T$:
\begin{eqnarray*}
T_{L} & = & \{ g \in T \mid g \cdot (x,y) = (x,y) ~if~ x \geq 0 \} \\ 
T_{R} & = & \{ g \in T \mid g \cdot (x,y) = (x,y) ~if~ x \leq 0 \} \\ 
T_{U} & = & \{ g \in T \mid g \cdot (x,y) = (x,y) ~if~ y \leq 0 \} \\ 
T_{D} & = & \{ g \in T \mid g \cdot (x,y) = (x,y) ~if~ y \geq 0 \}    
\end{eqnarray*}
Each of these groups is isomorphic to $F$, and so can be generated by two elements, which are necessarily small
as elements of $T$. 
Moreover $T = \langle T_{L}, T_{R}, T_{U}, T_{D} \rangle$.  It follows that $T$ is generated by $8$ elements,
each of which is small, such that the product of any two generators is also small.  This implies that $T$ fixes a point, 
by Lemma \ref{trees}, (4).
\end{proof}

\begin{figure} [!h]
\begin{center}
\includegraphics{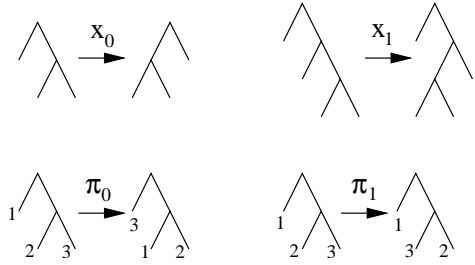}
\caption{A generating set for Thompson's group $V$.  The elements $x_{0}$, $x_{1}$, and $\pi_{0}$ generate Thompson's
group $T$.}
\label{myfig}
\end{center}
\end{figure}

\begin{theorem}
The group $V$ has Serre's property FA.
\end{theorem}

\begin{proof}
We recall from \cite{Cannon} that $T = \langle x_{0}, x_{1}, \pi_{0} \rangle$ and $V = \langle T, \pi_{1} \rangle$, where
$x_{0}$, $x_{1}$, $\pi_{0}$, and $\pi_{1}$ appear in Figure \ref{myfig}.  

Let $V$ act on a tree $\Gamma$.  By the previous theorem, we know that $T$ has a fixed point.  Thus, the elements $x_{0}$, $x_{1}$,
and $\pi_{0}$ all have fixed points, and any two of these elements will have a common fixed point.  The element $\pi_{1}$ has finite
order, so it must be elliptic.  It will be sufficient (by Lemma \ref{trees} (4)) 
to show that each product $\pi_{1} \pi_{0}$, $\pi_{1} x_{0}$, and $\pi_{1}x_{1}$ is elliptic.

First, we note that $\pi_{1} \pi_{0}$ is an element of finite order, so it must be elliptic.  
Next,
we note that $\pi_{1} x_{1}$ is small, so it must be elliptic.  

It is routine to check that $x_{0}^{-1} \pi_{1} x_{0} \pi_{1}$ is small, and therefore elliptic.  It follows that
$Fix( x_{0}^{-1}\pi_{1}x_{0} ) \cap Fix (\pi_{1}) \neq \emptyset$.  Now 
\begin{eqnarray*}
Fix(x_{0}^{-1}\pi_{1} x_{0} ) \cap Fix (\pi_{1}) \neq \emptyset & \Rightarrow & x_{0}^{-1} (Fix (\pi_{1})) \cap Fix (\pi_{1}) \neq \emptyset \\
& \Rightarrow & Fix(x_{0}^{-1}) \cap Fix (\pi_{1}) \neq \emptyset.
\end{eqnarray*}
This implies that $\pi_{1} x_{0}$ is elliptic.
\end{proof}

\bibliography{ends}

\begin{thebibliography}{10}

\bibitem{Brown}
Kenneth~S. Brown.
\newblock The homology of {R}ichard {T}hompson's group {$F$}.
\newblock In {\em Topological and asymptotic aspects of group theory}, volume
  394 of {\em Contemp. Math.}, pages 47--59. Amer. Math. Soc., Providence, RI,
  2006.

\bibitem{Cannon}
J.~W. Cannon, W.~J. Floyd, and W.~R. Parry.
\newblock Introductory notes on {R}ichard {T}hompson's groups.
\newblock {\em Enseign. Math. (2)}, 42(3-4):215--256, 1996.

\bibitem{iso}
Daniel~S. Farley.
\newblock Proper isometric actions of {T}hompson's groups on {H}ilbert space.
\newblock {\em Int. Math. Res. Not.}, (45):2409--2414, 2003.

\bibitem{actions}
Daniel~S. Farley.
\newblock Actions of picture groups on {CAT}(0) cubical complexes.
\newblock {\em Geom. Dedicata}, 110:221--242, 2005.

\bibitem{ross}
Ross Geoghegan.
\newblock {\em Topological Methods in Group Theory}.
\newblock to appear in the Springer Verlag series Graduate Texts in
  Mathematics.

\bibitem{GS}
Victor Guba and Mark Sapir.
\newblock Diagram groups.
\newblock {\em Mem. Amer. Math. Soc.}, 130(620):viii+117, 1997.

\bibitem{Higman}
Graham Higman.
\newblock {\em Finitely presented infinite simple groups}.
\newblock Department of Pure Mathematics, Department of Mathematics, I.A.S.
  Australian National University, Canberra, 1974.
\newblock Notes on Pure Mathematics, No. 8 (1974).

\bibitem{Hopf}
Heinz Hopf.
\newblock Enden offener {R}\"aume und unendliche diskontinuierliche {G}ruppen.
\newblock {\em Comment. Math. Helv.}, 16:81--100, 1944.

\bibitem{Houghton}
C.~H. Houghton.
\newblock Ends of locally compact groups and their coset spaces.
\newblock {\em J. Austral. Math. Soc.}, 17:274--284, 1974.
\newblock Collection of articles dedicated to the memory of Hanna Neumann, VII.

\bibitem{Klein}
Tom Klein.
\newblock {Filtered ends of infinite covers and groups}.

\bibitem{KR}
P.~H. Kropholler and M.~A. Roller.
\newblock Relative ends and duality groups.
\newblock {\em J. Pure Appl. Algebra}, 61(2):197--210, 1989.

\bibitem{NR2}
Terrence Napier and Mohan Ramachandran.
\newblock {Filtered ends, proper holomorphic mappings of Kahler manifolds to
  Riemann surfaces, and Kahler groups}.

\bibitem{NR1}
Terrence Napier and Mohan Ramachandran.
\newblock {Thompson's Group F is Not Kahler}.

\bibitem{Sageev}
Michah Sageev.
\newblock Ends of group pairs and non-positively curved cube complexes.
\newblock {\em Proc. London Math. Soc. (3)}, 71(3):585--617, 1995.

\bibitem{Serre}
Jean-Pierre Serre.
\newblock {\em Trees}.
\newblock Springer Monographs in Mathematics. Springer-Verlag, Berlin, 2003.
\newblock Translated from the French original by John Stillwell, Corrected 2nd
  printing of the 1980 English translation.

\end{thebibliography}
\nocite{*}
\bibliographystyle{plain}

\end{document}